\documentclass[11pt]{article}
\usepackage{amsmath,amssymb}
\newtheorem{theorem}{Theorem}
\newtheorem{lemma}{Lemma}
\newtheorem{proposition}{Proposition}

\newtheorem{remark}{Remark}
\def\s{\sigma}
\newcommand{\baa}{\begin{eqnarray}}
\newcommand{\eaa}{\end{eqnarray}}
\newcommand{\baaa}{\begin{eqnarray*}}
\newcommand{\eaaa}{\end{eqnarray*}}

\newcommand {\Var}{\mbox{Var}}

\def\defi{\stackrel{{\scriptscriptstyle \Delta}}{=}}

\def\a{\alpha}

\def\o{\omega}
\def\O{\Omega}

\def\F{{\cal F}}

\def\Ind{{\,\rm Ind\,}}
\def\Ind{{\mathbb{I}}}

\def\esssup{\mathop{\rm ess\, sup}}

\def\Var{{\rm Var\,}}

\def\R{{\bf R}}
\def\E{{\bf E}}
\def\P{{\bf P}}

\def\b{\beta}
\def\s{\delta}

\def\oo{\bar}
\def\s{\sigma}

\def\oo{\bar}

\newtheorem{condition}{Condition}

\title{
On limit  periodicity of discrete time stochastic processes
}
\begin{document}
\author{Alexandra Rodkina
\footnote{ Department of Mathematics,  The University of the West
Indies, Kingston 7, Jamaica }, Nikolai Dokuchaev\footnote{Department
of Mathematics and Statistics, Curtin University, Australia}, and
John Appleby\footnote{School of Mathematical Sciences, Dublin City
University,  Dublin 9, Ireland}\footnote{ The first and the second
author were supported  by ARC grant of Australia DP120100928.}}

\maketitle

\begin{abstract}
We consider a discrete time dynamic system described by a difference equation
with periodic  coefficients and with additive stochastic noise.
 We investigate the possibility of the  periodicity for the solution. In particular, we
 found sufficient conditions for existence of a
periodic process such that the solution converges to it, including
almost surely convergence.
\\
 {\bf Key words}: discrete time dynamic systems, stochastic difference equations, periodic
 solutions.
 \\
{\it AMS  classification:} 
  39A23,  
        39A30,  
    39A50    
\end{abstract}

\section{Introduction and problem setting}
The a.s. asymptotic stability of stochastic discrete time processes
has been widely addressed; see, e.g., \cite{AMaR}\index{, \cite{AR},
\cite{ABR2}, \cite{ABR1},
 \cite{B},}-
\cite{Ma}, \cite{O} and the bibliography here). However, periodicity
of these processes are not so well represented. In
\cite{M},\cite{U}, conditions of  periodicity in the distributions
were obtained for discrete time systems; a review of periodicity for
nonlinear discrete time equations can be found in \cite{E}. In
\cite{I}, \cite{Mar},\cite{2}, conditions of periodicity in the
distributions were obtained for continuous time systems. In this
article we obtain sufficient conditions for convergence of solutions
 to a periodic process  in a strong sense. We
consider processes  described by stochastic difference equations
with periodic  coefficients and with decaying additive stochastic
noise.
 We investigate  the asymptotic properties.
In particular, we found some sufficient conditions  when this
convergence can be achieved almost surely.

   Let $(\Omega, {\mathcal{F}},\P)$ be a complete  probability space, with elementary events  $\o\in \O$.

 Consider a stochastic process $X_n=X_n(\o)$,  $\o\in \O$,  that evolves as
  \baaa
&&X_{n+1}=X_n+a_nX_n+\sigma_n \xi_{n+1}, \quad n=0,1,2,.., \quad
\\&&X(0)=X_0. \label{eq:intr} \eaaa Here $X_0\in\R$,
$\{\sigma_n\}$ is a nonrandom sequence of real
numbers, $\{a_n\}=\{a_n(\o)\}$ is a random  periodic sequence, and $\{\xi_n\}=\{\xi_n(\o)\}$ is a sequence of random variables.

 We
assume that $\{a_n\}$ is a
 periodic sequence of bounded random  variables  with a period $K$, i.e. $a_{n+K}=a_n$ a.s. for all $n=0,,1,2,...$.

We assume that  $\sigma_n\to 0$ as $n\to +\infty$.

 We assume that there exists
$p\in[1,+\infty]$ such that $\|\xi_k\|_{L_p(\O)}\le 1$ for all $k$.

For the brevity, we denote $L_p(\O)=  L_p(\Omega, \mathcal F, \P)$.

We use the standard abbreviation ``{\it a.s.}'' for the wordings
``almost sure'' or ``almost surely'' throughout the text.

\section{The main results}
\def\oo{\overline}
Let $L\defi (1+a_0)(1+a_1)\cdots (1+a_{K-1})$.
\begin{theorem}  \label{ThM}
\begin{enumerate}
\item
 If $\esssup_{\o}|L|<1$ then $\|X_n\|_{L_p(\O)}\to 0$ as
$n\to +\infty$. \item Assume that $p\ge 2$, $\esssup_{\o}|L|>1$,   $\xi_n$ are
mutually independent and also independent on $\{a_n\}$, and either
$X_0\neq 0$ or $\inf_n\Var\xi_n>0$. Then ${\rm limsup}
\|X_n\|_{L_p(\O)}\to +\infty$ as $n\to +\infty$.
\item
If $L=1$ a.s. then  $\|X_n-X_{n+K}\|_{L_2(\O)}\to  0$ as $n\to +\infty$.
 \item
If $L=1$ a.s. and $\lim_{n\to \infty}\sigma_n\xi_{n+1}=0$ a.s. then
$X_n-X_{n+K}\to  0$ a.s. as $n\to +\infty$.
\end{enumerate}
\end{theorem}
\begin{remark}
\label{rem:sigmaxi0} Condition when,  a.s., $\lim_{n\to
\infty}\sigma_n\xi_{n+1}=0$, are given e.g. in \cite{ABR1, ABR2}.
\end{remark}

Starting from now and up to the end of this paper, we assume that
$L=1$ and that the following conditions is satisfied:
\begin{condition}
\label{con0}  One of the following conditions is satisfied:
\begin{enumerate}
\item
 $\sum_{k=0}^\infty |\s_k|<+\infty$; or
\item $p=2$,
 $\sum_{k=0}^\infty |\s_k|^2<+\infty$, and $\{\xi_n\}$ is a sequences of  independent on $\{a_n\}$ and mutually independent
identically distributed random variables such that $\E\xi_n=0$.
\end{enumerate}
\end{condition}

\begin{lemma}
\label{lemma0} The sum \baaa \oo Y=\sum_{k=0}^{\infty}\s_k\xi_k.
\eaaa belongs to  $\oo Y\in L_p(\O,\F,\P)$; it is defined as the
limit of the partial sums in this space.
\end{lemma}

Let $b_{k,n}=\prod_{i=k}^{n-1}(1+a_i)$, $k<n$, $b_{n,n}=1$, and  let \baaa \label{def:J}  \bar Q_n=b_{0,n}\bar Y. \eaaa Note that
$\bar Q_n$ is a.s. a $K$-periodic process.

\begin{theorem}
\label{ThM1}
 $\lim_{n\to 0}\|X_n-\bar Q_n\|_{L_p(\O)}=0$.
\end{theorem}

The  following two question arises: is it true that $\tilde X_n$
such that   $X_n-\oo Q_{n}\to  0$ as $n\to +\infty$ a.s.? This
question is addressed in the following theorem.
\begin{theorem}
\label{thm:period}  Let  at least one  of the following conditions
is satisfied:
\begin{enumerate}
\item
$\sup_{n\ge 0,\o\in\O}|\xi_n(\o)|<+\infty$,  or
\item Condition \ref{con0}(ii) is satisfied.
\end{enumerate}
Then
 $\lim_{n\to 0}|X_n-\bar Q_n|=0$ a.s. and $\lim_{n\to 0}|X_n-X_{n+K}|=0$ a.s.
\end{theorem}

Note that the assumption (i) in Theorem \ref{thm:period} does not require that $\xi_n$ are independent and independent from $\{a_k\}$.

\section{Proofs}
 Let
 \baaa
\psi_{n,m}=b_{n+1,n+m}\s_{n+1}\xi_{n+1}+b_{n+2,n+m}\s_{n+2}\xi_{n+2}+...+b_{n+m,n+m}\s_{n+m}\xi_{n+m}.
\eaaa
 {\em Proof of Theorem
\ref{ThM}}.  We have that \baaa
X_n=b_{0,n}X_0+b_{1,n}\s_1\xi_1+b_{2,n}\s_2\xi_2+...+b_{n,n}\s_n\xi_n.
\eaaa Let $a_k=a_{k+K}$, $b_{0,K}=L$ for some $K$.  In this case,
$b_{k,K+k}=L$ for all $k$.

We have that \baaa
X_{n+K}=b_{0,n+K}X_0+b_{1,n+K}\s_1\xi_1+b_{2,n+K}\s_2\xi_2+...+b_{n+K,n+K}\s_{n+K}\xi_{n+K}
\\=b_{0,n+K}X_0+b_{1,n+K}\s_1\xi_1+b_{2,n+K}\s_2\xi_2+...+b_{n,n+K}\s_{n}\xi_{n}
+\psi_{n,K}
\\=L(b_{0,n}X_0+b_{1,n}\s_1\xi_1+b_{2,n}\s_2\xi_2+...+b_{n,n}\s_{n}\xi_{n})
+\psi_{n,K}. \eaaa It follows that \baa X_{n+K}=LX_n +\psi_{n,K}.
\label{XnK}\eaa

 Let $Y_i=X_{iK}$ and $\eta_i=\psi_{iK,K}$. By (\ref{XnK}), it follows that \baa Y_{i+1}=LY_{i}+\eta_i,\quad
i\ge 0,\qquad  Y_0=X_0. \label{Yi}\eaa

Note that the set $\{b_{n+s,n+K}, s\in\{1,...,K\}, n>0 \}$ is
bounded in $L_\infty(\O)$. By the assumptions,   $\s_k\to 0$ as
$k\to +\infty$.  Hence  $\|\eta_i\|_{L_2(\O)}\to 0$ as $i\to
+\infty$ and $\|\psi_{n,K}\|_{L_2(\O)}\to 0$   as $n\to +\infty$.
Then statement (iii) follows.  Under the conditions of
statement (iv), $\psi_{n,K}\to 0$ a.s.   as $n\to +\infty$. Then
statement  (iv) follows.

Let us prove statement (i). Let $\a_i=\|Y_{i+1}\|_{L_2(\O)}$,
$\b_i=\|\eta_i\|_{L_2(\O)}$, $L_0=\esssup_{\o}|L|$. We have that
\baaa \a_{i+1}\le L_0\a_{i}+\b_i,\quad i\ge 0,\qquad  \a_0=|X_0|.
\label{Yii}\eaaa Consider the equation \baaa \oo\a_{i+1}=
L_0\oo\a_{i}+\b_i,\quad i\ge 0,\qquad  \oo\a_0=|X_0|.
\label{Yiii}\eaaa By the properties of the solutions of this
equation, we have that $\oo a_i\to 0$ as $i\to +\infty$. Since $0\le
\a_i\le \oo\a_i$, we obtain  statement (i).

Let us prove statement (ii). Consider event $A=\{|L|>1\}$. By the assumptions, $\P(A)>0$.   We  have that
\baaa
Y_{i+1}=L^i\eta_0+L^{i-1}\eta_1+\cdots + L\eta_{i-1}+\eta_i
\eaaa
and
\baaa
\Ind_AY_{i+1}=\Ind_AL^i\eta_0+\Ind_AL^{i-1}\eta_1+\cdots \Ind_AL\eta_{i-1}+\Ind_A\eta_i.
\eaaa
Note that $\Ind_A$ is non-random on the conditional probability space given $\{a_k\}$, and
the random variables $\Ind_AL^{m-1}\eta_m$ are independent on the conditional probability space given $\{a_k\}$.
Hence
\baaa
&&\Var(\Ind_AY_{i+1}|(a_k))\\&&=\Var(\Ind_AL^i\eta_0|(a_k))+\Var(\Ind_AL^{i-1}\eta_1|(a_k))+\cdots \Var(\Ind_AL\eta_{i-1}|(a_k))+\Var(\Ind_A\eta_i|(a_k))
\\
&&=L^i\Var(\Ind_A\eta_0|(a_k))+L^{i-1}\Var(\Ind_A\eta_1|(a_k))+\cdots L\Var(\Ind_A\eta_{i-1}|(a_k))+\Var(\Ind_A\eta_i|(a_k))
\\
&&=\Ind_A(L^i\Var(\eta_0|(a_k))+L^{i-1}\Var(\eta_1|(a_k))+\cdots L\Var(\eta_{i-1}|(a_k))+\Var(\eta_i|(a_k)).
\eaaa
Here $\Var(\cdot |(a_k))$ is
the conditional variance  given $\{a_k\}$. By the assumptions, there exists $m$ such that $\Var \xi_m>0$.  If $|L|>1$ then all $b_{k,n}\neq 0$. Hence there exists $j$ such that
 $\Var (\eta_j|(a_k))>0$ a.s. given that $|L|>1$.
 We obtain immediately that
$\Var (Y_i|(a_k))\to +\infty$  as $i\to +\infty$ a.s. given that $|L|>1$  Hence $\E (Y_i^2|(a_k))\to +\infty$ as $i\to +\infty$ a.s.. (This is the conditional second moment under the same condition).   Therefore,  $\E Y_i^2\to +\infty$ as $i\to +\infty$.    This completes the
proof of Theorem \ref{ThM}.
\index{
Let us prove statement (ii).  We obtain immediately that the
conditional variance $\Var(X_n|(a_k), |L|>1)\to 0$ as $n\to +\infty$. Here $\Var(X_n|(a_k), |L|>1)$ is
conditional variance  given $(a_k)$ and given that $|L|>1$. Hence $\E(X_n^2|(a_k), |L|>1)\to +\infty$ as $n\to +\infty$  a.s.
(This is the conditional second moment under the same condition).    This completes the
proof of Theorem \ref{ThM}.} $\Box$

\par {\em Proof of Lemma \ref{lemma0}}. Let \baaa\oo Y_{i}\defi \sum_{j=0}^{i}\s_j\xi_j.\eaaa
We have that
\baaa
\oo Y_{i}-\oo Y_{i+m}=  \sum_{j=i+1}^{i+m}\s_j\xi_j.
\eaaa
  Let us assume first that Condition \ref{con0}(i) is satisfied.
\baaa
\|\oo Y_{i}-\oo Y_{i+m}\|_{L_p(\O)}=  \Bigl\|\sum_{j=i+1}^{i+m}\s_j\xi_j\Bigr\|_{L_p(\O)}\le   \sum_{j=i+1}^{i+m}|\s_j|\|\xi_j\|_{L_p(\O)}\le
\sup_{k\ge 0}\|\xi_k\|_{L_p(\O)}\sum_{j=i+1}^{i+m}|\s_j|\\\le \sum_{j=i+1}^{\infty}|\s_j|.
\eaaa
By
Condition \ref{con0}(i), $\sum_{j=i+1}^{\infty}|\s_j|\to 0$ as $j\to +\infty$. Hence $\{\oo Y_{i}\}$ is a Cauchy sequence in $L_p(\O)$.
Then the statement of lemma follows  in this case.

Let us assume first that Condition \ref{con0}(ii) is satisfied.
Let \baaa\oo Y_{i}\defi \sum_{j=0}^{i}\s_j\xi_j.\eaaa
We have that $\E Y_i=0$ and
\baaa
\oo Y_{i}-\oo Y_{i+m}=  \sum_{j=i+1}^{i+m}\s_j\xi_j.
\eaaa
We have that
\baaa
\E(\oo Y_{i}-\oo Y_{i+m})^2=\Var(\oo Y_{i}-\oo Y_{i+m})=  \Var\Bigl(\sum_{j=i+1}^{i+m}\s_j\xi_j\Bigr)=   \sum_{j=i+1}^{i+m}|\s_j|^2\Var(\xi_j)\\\le
\sup_{k\ge 0}\|\xi_k\|_{L_2(\O)}\sum_{j=i+1}^{i+m}|\s_j|^2\le \sum_{j=i+1}^{\infty}|\s_j|^2.
\eaaa
By
Condition \ref{con0}(ii), $\sum_{j=i+1}^{\infty}|\s_j|^2\to 0$ as $j\to +\infty$. Hence $\{\oo Y_{i}\}$ is a Cauchy sequence in $L_2(\O)$.
Then the statement of lemma follows  in this case.
This completes the
proof of Lemma \ref{lemma0}.  $\Box$

{\em Proof of Theorem \ref{ThM1}}.  Let an integer $m\in (0,K]$ be
given. We have that
$$
X_{nK+m}=b_{nK,m}X_{nK}+\psi_{nK,m}=b_{0,m}X_{nK}+\psi_{nK,m}.
$$
 It gives
 \baa
X_{nK+m}=b_{0,m}X_{nK}+\psi_{nK,m}=
b_{0,m}Y_{n}+\psi_{nK,m}\nonumber\\=b_{0,m}\bar Y +\psi_{nK,m}+
b_{0,m}(Y_{n}-\bar Y). \label{Xk}\eaa Clearly, $\|b_{0,m}(Y_{n}-\bar
Y)\|_{L_p(\O)}\to 0$. By the definition, $Q_{nK+m}=b_{0,m}\bar Y$.
Further, we have   that
 $\|\psi_{nK,m}\|_{L_p(\O)}\to 0$.
This completes the proof of Theorem \ref{ThM1}. $\Box$

Up to the end of this paper, we assume that the assumptions of
Theorem \ref{thm:period} are satisfied.

\def\ssigma{\bar\sigma}

\begin{lemma}
\label{lem:limdistr} Assume that condition \ref{con0}(ii) is satisfied. In this case,
 \begin{enumerate}
 \item [(i)] there exists a.s. finite random variable $\bar Y$ such that, a.s.,
 \baaa
\label{def:barX} \bar Y=\lim_{i\to \infty}Y_i. \eaaa

\item[(ii)]  $\lim_{i\to \infty}\mathbf E |Y_i-\bar Y|^2=0$ and $\mathbf E |\bar Y|^2<\infty$.

\item[(iii)]
 $
Var \bar Y=\lim_{n\to \infty}\langle
Y_n\rangle=\sum_{i=1}^\infty\sigma_i^2<\infty. $
\end{enumerate}
\end{lemma}

We define by $\ell_2$ a Banach space of sequences  $\{y_n\}$ of real numbers, such that
\[
\|y\|_{\ell_2}=\sum_{n=0}^\infty y_n^2<\infty.
\]
  If $\{M_n\}_{n\in \mathbf N}$ is a square integrable martingale,  $M_0=0$ and
$M_n=\sum_{i=1}^n \rho_i$, where  $\{\rho_n\}_{n \in \mathbf N}$ is
a sequence of independent random variables with $\mathbf E\rho_n=0$
and $\mathbf E\rho^2_n<\infty$, then the \emph{quadratic variation}
$\langle M_n\rangle$    of $M$ is  defined by
\[
\langle M_n\rangle=\mathbf E M^2_n=\sum_{i=1}^n\mathbf{E}\rho_i^2.
\]
In this situation the process $\langle M_n\rangle$ is not random and
coincides with the variance of $M_n$.

 A detailed exposition of the
definitions and facts of the theory of random processes can be found
in, for example, \cite{Shiryaev96}.

The Proposition below is a variant of the martingale convergence
theorem (see e.g. \cite{Shiryaev96}).
\begin{proposition}
\label{pr:limdistr} Let condition \ref{con0}(ii) be satisfied, let $\{\ssigma\}\in \ell_2$, and let
 \begin{equation}
  \label{def:martM}
  M_n=\sum_{i=0}^{n-1} \ssigma_i\xi_{i+1}.
  \end{equation}
Then
 \begin{enumerate}
 \item [(i)] there exists a.s. finite random variable $\bar M$ such that, a.s.,
 \baaa
\label{def:barM} \bar M=\lim_{i\to \infty}M_i. \eaaa

\item[(ii)]  $\lim_{i\to \infty}\mathbf E |M_i-\bar M|^2=0$ and $\mathbf E |\bar M|^2<\infty$.

\item[(iii)]
 $
Var \bar M=\lim_{n\to \infty}\langle M_n\rangle=\sum_{i=1}^\infty\ssigma_i^2<\infty.
$
\end{enumerate}
\end{proposition}

{\em Proof of Proposition \ref{pr:limdistr}}.
Item  (i) is a variant of martingale convergence theorems
(see e.g. \cite{Shiryaev96}).

For $n>k$, we have \baaa \|M_n-M_k\|_{ L_2(\O)}= &\mathbf E
\left|\sum_{i=k}^{n-1}\ssigma_i \xi_{i+1} \right|^2\le
\sum_{i=k}^{n-1}\ssigma_i^2\to 0, \quad \hbox{as} \quad k\to \infty,
\eaaa It follows that $\{M_n\}_{n\in \mathbf N}$ is a Cauchy
sequence in $ L_2(\O)$. Fix $\varepsilon>0$ and choose $N\in \mathbf
N$ such that
\[
\|M_n-M_k\|_{ L_2(\O)}<\varepsilon, \quad \hbox{as} \quad n, k\ge N.
\]
Then, by the Fatou's lemma and the fact that $M_k\to \bar M$ a.s.,
we have, for $n\ge N$, \baaa &&\|M_n-\bar M\|_{ L_2(\O)}=\mathbf
E|M_n-\bar M|^2=\mathbf E|M_n-\lim_{k\to \infty}M_k|^2\\&&= \mathbf
E \left\{\lim_{k\to \infty} |M_n-M_k| ^2 \right\}=\mathbf E
\left\{\liminf_{k\to \infty} |M_n-M_k| ^2 \right\}\\&&\le
\liminf_{k\to \infty}\mathbf E \left\{ |M_n-M_k| ^2
\right\}=\liminf_{k\to \infty}\|M_n-M_k\|_{ L_2(\O)} <\varepsilon.
\eaaa By Minkosvki inequality this also implies that $\bar M\in
L_2(\O)$, which completed the proof of (ii).

To prove (iii) we estimate, \baaa &&|\langle M_n\rangle-Var \bar M|
=|\mathbf EM_n^2-\mathbf E\bar M^2|=|\mathbf E(M_n^2-\bar
M^2)|\\&&=|\mathbf E\left[(M_n-\bar M) (M_n+\bar M)\right]|\le
 \sqrt{\mathbf E|M_n-\bar M|^2}\sqrt{\mathbf E|M_n+\bar M|^2}\\&&\le K\sqrt{\mathbf E|M_n-\bar M|^2}\to 0 \quad \hbox{as $n\to \infty$}.
\eaaa This can be proved differently. We estimate \baaa &&\|M_n\|_{
L_2(\O)} \le \|\bar M\|_{ L_2(\O)} +\|\bar M-M_n\|_{ L_2(\O)},\\
&&\|\bar M\|_{ L_2(\O)} \le \|M_n\|_{ L_2(\O)} +\|\bar M-M_n\|_{
L_2(\O)},
 \eaaa
 which implies that
 \[
-\|\bar M-M_n\|_{ L_2(\O)}\le \|M_n\|_{ L_2(\O)} - \|\bar M\|_{
L_2(\O)} \le\|\bar M-M_n\|_{ L_2(\O)}.
  \]
 So, as $n\to \infty$,
 \[
 |Var \bar M-\langle M_n\rangle|=\left| \|\bar M\|_{ L_2(\O)}-\|M_n\|_{ L_2(\O)}\right|\le \|\bar M-M_n\|_{ L_2(\O)}\to 0.
  \]

 \par
{\em Proof  of Lemma \ref{lem:limdistr}}.
 By (\ref{Yi}),  we have that $Y_{i+1}=X_0+\sum_{j=0}^i\eta_i=X_0+\sum_{n=0}^{(i+1)K}\s_n\xi_{n}$. Hence
$\{Y_i\}_{i\ge 0}$ is a subsequence of $\{X_0+M_n\}_{n\ge 0}$. Then the  proof follows immediately from Proposition \ref{pr:limdistr} applied with $\ssigma_n=\sigma_n$.

{\em Proof  of Theorem \ref{thm:period}}. Let $m\in(0,K]$. (i) Let
assumption (i) holds.  In this case,   \baaa |Y_{n}-\bar Y|\le
\sup_{n,\o}|\xi_n|\sum_{k=n}^\infty|\s_k|\to 0\quad\hbox{as}\quad n\to
+\infty \quad\hbox{a.s.},\eaaa since  $\sum_{k=0}^\infty
|\s_k|<+\infty$, by Condition \ref{con0}. It follows that $b_{0,m}(Y_{n}-\bar Y)\to 0$ a.s. By
the definition, $Q_n=b_{1,m}\bar Y$. Similarly, \baaa
|\psi_{nK,m}|\le
K\sup_{n,\o}(1+|a_k|)^K\sup_{n,\o}|\xi_n|\sum_{k=n}^\infty|\s_k|\to
0\quad\hbox{as}\quad n\to +\infty\quad\hbox{a.s.} \eaaa By (\ref{Xk}), the proof of Theorem
\ref{thm:period} follows for this case.
\par
Let assumption (ii) holds.
 By Lemma \ref{lem:limdistr}, it follows that
$b_{0,m}(Y_{n}-\bar Y)\to 0$ a.s. By the definition,
$Q_n=b_{1,m}\bar Y$. Further, we obtain that
 $\psi_{nK,m}\to 0$ a.s.,  by Proposition \ref{pr:limdistr}
applied with
$\ssigma_{n+k}=b_{n+K,n+m}\s_{n+k}$  on the conditional probability space given   $\{a_k\}$.  By (\ref{Xk}), the proof of Theorem \ref{thm:period} follows. $\Box$

\begin{remark} The statements
and the proofs of Theorems \ref{ThM} - \ref{thm:period} hold without changes for the case where $X_n$ and
$\xi_n$ are $d$-dimensional processes for some $d>1$, where $a_n$ and
$\s_n$ are $d\times d$-dimensional matrices, and where
$(I+a_{K-1})(I+a_{K-2})\cdots (I+a_0)=LI$. Here $I$ is the unit
matrix in $R^{d\times d}$, $L\in \R$.
\end{remark}
\index{
\begin{remark} As can be seen from the proofs, the result of
 Theorem \ref{thm:period} is based on the a.s. convergence of
 $\sum_{n=0}^{\infty}\ssigma_n\xi_{n}$ established in Proposition \ref{pr:limdistr}. So far, we were unable to establish if this convergence can be achieved
 for  the case where
 $\sum_{n=0}^{\infty}|\s_n|^2=+\infty$. We leave it for  future research.
\end{remark}
}

\end{document}